# EXPLICIT EXPONENTIAL MAPS FOR HECKE CHARACTERS AT ORDINARY PRIMES

LI GUO

ABSTRACT. Let $E$ be an elliptic curve with complex multiplication by the ring of integers of an imaginary quadratic field $K$. The theory of complex multiplication associates $E$ with a Hecke character $\psi$. The Hasse-Weil $L$-function of $E$ equals the Hecke $L$-function of $\psi$, whose special value at $s = 1$ encodes important arithmetic information of $E$, as predicted by the Birch and Swinnerton-Dyer conjecture and verified by Rubin[**Ru**] when the special value is non-zero. For integers $k, j$, special values of the Hecke $L$-function associated to the Hecke character $\psi^k \bar\psi^j$ should encode arithmetic information of the Hecke character $\psi^k \bar\psi^j$, as predicted by the Bloch-Kato conjecture[**B-K**]. When $p$ is a prime where $E$ has good, ordinary reduction, the $p$-part of the conjecture has been verified when $j = 0$[**Ha**] and when $j \neq 0, p > k$[**Gu**]. To verify the conjecture in other cases, it is important to have an explicit description of the exponential map of Bloch and Kato. In this paper we provide such an explicit exponential map for the case $j \neq 0$, $p \nmid k$.

Let $K$ be an imaginary quadratic field with ring of integers $\mathcal{O}_K$. Let $E$ be an elliptic curve with complex multiplication by $\mathcal{O}_K$, and let $\psi$ be the Grössencharacter associated to $E$ by the theory of complex multiplication. For integers $k, j$, we have the Grössencharacter $\varphi = \psi^k \bar\psi^j$. The behavior of the Hecke $L$-function $L(\varphi^{-1}, s)$ at $s = 0$ is predicted by the conjecture of Bloch and Kato. When $k > -j \geq 0$, $L(\varphi^{-1}, 0) \neq 0$ and $E$ has good, ordinary reduction at $p$, the $p$-part of the Bloch-Kato conjecture is verified by using Iwasawa theory when $j = 0$[**Ha**] and when $j \neq 0, p > k$[**Gu**]. In order to study the remaining case (when $j \neq 0$, $p \leq k$), it is important to get an explicit exponential map for the $p$-adic representation from the Grössencharacter $\varphi$. Such an explicit exponential map is needed to describe the image of the exponential map and thus to get control of the Tamagawa number that appears in the Bloch-Kato conjecture. We will provide such an explicit exponential map (Theorem 2) in this paper.

When $j = 0$, an explicit exponential map is given by Harrison[**Ha**], building on previous work of Bloch and Kato [**B-K**], and Fontaine [**Fo2**]. It is obtained by first finding an explicit reciprocity law for the height one Lubin-Tate formal groups and then using the fact that the $p$-adic representation from the Grössencharacter actually comes naturally from such a Lubin-Tate formal group. This is not the case when $j \neq 0$. However we will show that for $p \nmid k$, it is still possible to relate the $p$-adic character from $\varphi$ to a formal group over a finite extension of $K$. Then an explicit exponential map can be obtained (Theorem 2). Once this is done, the $p$-part of the Bloch-Kato conjecture could be proved by an argument that is lengthy

---

*Date*: July 25, 1997.
1991 *Mathematics Subject Classification.* Primary 11S31; Secondary 11G45, 11R23, 14G10.
This research is supported in part by an NSF Grant.

and technical, but otherwise analogous to earlier work[**Ha**, **Gu**]. Thus we have chosen not to give the full details here.

For good, ordinary primes $p$, the $p$-part of the Bloch-Kato conjecture of $\varphi$ has been verified in [**Gu**] when $p > k$. Thus the explicit exponential map developed here is essential only for the finite many primes $p \leq k, p \nmid k$. The same method could be applied to the case when the character $\varphi$ is twisted by a finite character, i.e., to the case when we have a Grössencharacter over an imaginary quadratic field of class number one. More interesting, of course, is the case when the prime $p$ is supersingular. This time the height of the formal group is two and until recently the only known explicit reciprocity law had been the one of Wiles[**Wi**] when $k = 1, j = 0$. Rubin[**Ru**] used this explicit reciprocity law and the Iwasawa main conjecture to verify the Birch and Swinnerton-Dyer conjecture in an important case. It is equivalent to the Bloch-Kato conjecture for the character $\psi$ when $L(\psi^{-1}, 0) \neq 0$. In an important paper[**Ka**], Kato formulated a very general explicit reciprocity law for Lubin-Tate formal group of any height. It apparently implies the explicit exponential map for the $p$-adic characters from $\psi^k$, $p$ supersingular. It is our hope that, combining with Kato's result, the method of this paper would shed further light on the study of other characters at supersingular primes.

In §1, notations will be introduced and the explicit reciprocity law of Harrison will be stated, together with other preliminary results. In §2 we prove a relation between the p-adic Hecke characters and the formal group characters for Lubin-Tate groups of height one (Proposition 2.1). With this relation we obtain an explicit exponential map for the $p$-adic characters (Theorem 2), which, in turn, enable us to, in essence, describe the image of the exponential map (Proposition 2.2).

It should be clear that this work is based the ideas from [**Fo2**, **Ha**]. The author is indebted to their work. The author is also grateful to Karl Rubin for helpful discussions.

## 1. Notations and preliminaries

For each prime $p$ of $\mathbb{Q}$, let $\bar{\mathbb{Q}}_p$ be the algebraic closure of $\mathbb{Q}_p$, and let $\mathbb{Q}_p^{\text{un}}$ be the maximal unramified extension of $\mathbb{Q}_p$ in $\bar{\mathbb{Q}}_p$. Denote $\mathbb{C}_p$ (resp. $\mathbb{D}_p$) for the completion of $\bar{\mathbb{Q}}_p$ (resp. $\mathbb{Q}_p^{\text{un}}$) and denote $\mathcal{O}_{\mathbb{C}_p}$ and $\mathcal{O}_{\mathbb{D}_p}$ for the corresponding rings of integers. Let $\pi$ be a uniformiser of $\mathbb{Q}_p$ and $\mathbf{F}$ be a one-dimensional Lubin-Tate module over $\mathbb{Z}_p$ of height one for $\pi$. Let $\mathfrak{m}_p$ be the maximal ideal of $\mathcal{O}_{\mathbf{C}_p}$ and let $\mathbf{F}_{\pi^n}$ be the group of $\pi^n$-torsion points of $\mathbf{F}(\mathfrak{m}_p)$. Let $L$ be the finite unramified extension of $\mathbb{Q}_p$ of degree $d$. Let $F$ be the Frobenius. Then $L(\mathbf{F}_{\pi^\infty})/L$ is a totally ramified $\mathbb{Z}_p^\times$-extension. Thus the action of $G_{\mathbb{Q}_p}$ on $\mathbf{F}_{\pi^\infty}$ gives rise to a continuous Galois character $\kappa : G_{\mathbb{Q}_p} \to \mathbb{Z}_p^\times$ with kernel $G_{L(\mathbf{F}_{\pi^\infty})}$. Let $\mathbb{Q}_p(\kappa^n)$ be the 1-dimensional $G_L$-representation over $\mathbb{Q}_p$ induced from $\kappa^n$, $n \in \mathbb{Z}$. Harrison[**Ha**] proved an explicit reciprocity law for $\mathbb{Q}_p(\kappa^n)$. We need more notation to formulate it. For more details see [**de**] and [**Ha**].

Let $\mathbb{G}_m$ be the formal multiplicative group. There is $\theta(T) \in \mathcal{O}_{\mathbf{D}_p}[[T]]$ giving a $\mathbb{Z}_p$-isomorphism from $\mathbb{G}_m$ onto $\mathbf{F}$ over $\mathcal{O}_{\mathbf{D}_p}$. We fix a coherent set of primitive $p$-power roots of unity $(\zeta_n)_{n\geq 0}$, that is, $\zeta_0 = 1, \zeta_1 \neq 1, \zeta_{n+1}^p = \zeta_n$ for $n \geq 0$. Define $\pi_n = \theta^{F^{-n}}(\zeta_n - 1), n \geq 0$. Then we have $\pi_n \in \mathbf{F}_{\pi^n}/\mathbf{F}_{\pi^{n-1}}$ and $[\pi]_\mathbf{F}(\pi_n) = \pi_{n-1}$. Let $\lambda_\mathbf{F}(T) \in \text{Hom}(\mathbf{F}, \mathbb{G}_a)$ be the formal logarithm of $\mathbf{F}$ and let $\log(1 + T) = \sum_{i=0}^\infty \frac{(-1)^{i-1}}{i} T^i$ be the formal logarithm of $\mathbb{G}_m$. Then $\lambda_\mathbf{F}(\theta(T)) = \Omega_p \log(1 + T)$.

Here $\Omega_p \in \mathcal{O}_{\mathbf{D}_p}$ satisfies $\Omega_p^\sigma = (\kappa(\sigma)/\chi_p(\sigma))\Omega_p$, $\sigma \in G_{\mathbb{Q}_p}$, and $\chi_p$ is the $p$-cyclotomic character associated to the formal group $\mathbb{G}_m$.

Now write $L_n$ for $L(\mathbf{F}_{\pi^n})$, $n > 0$, $L_0 = L$, and $U_{L_n}^1$ for the units of $\mathcal{O}_{L_n}$ which are congruent to 1 modulo the maximal ideal. Define $U_\infty^1 = \varprojlim U_{L_n}^1$ where the inverse limit is with respect to the norm maps. For a norm coherent sequence $u = (u_n) \in U_\infty^1$ of 1-units, let $g_u(X) \in \mathcal{O}_L[[X]]^\times$ be the Coleman power series such that $g_u^{F^{-i}}(\pi_i) = u_i$, $i \geq 1$ and $g_u(0)^{1-F^{-1}} = u_0$.

For $k \geq 1$, the $k$-th Coates-Wiles homomorphism $\phi_{CM}^k : U_\infty^1 \to L$ is defined by

$$\phi_{CW}^k(u) = \left[\frac{1}{\lambda'_{\mathbf{F}}(Y)} \frac{d}{dY}\right]^k \log g_u(Y)|_{Y=0}.$$

$\phi_{CW}^k$ is a $\mathbb{Z}_p$-homomorphism and $\phi_{CW}^k(u^\sigma) = \kappa^k(\sigma)\phi_{CW}^k(u)$ for $\sigma \in G_L$.

Next we define the reciprocity map. Define $t_\pi = \Omega_p t \in B_{\mathrm{crys}}^+$ where $t$ is the usual element of $B_{\mathrm{crys}}$ on which $G_{\mathbb{Q}_p}$ acts via the $p$-adic cyclotomic character. Then we have $t_\pi^F = \pi t_\pi$ and $t_\pi^\sigma = \kappa(\sigma) t_\pi$, $\sigma \in G_L$. Multiplying each term of the exact sequence

$$0 \to \mathbb{Q}_p(k) \to J_{\mathbb{Q}}^{[k]} \xrightarrow{1-p^{-k}F} B_{\mathrm{crys}}^+ \to 0$$

from [**B-K**, (1.13)] by $\Omega_p^k$ we get the exact sequence

$$0 \to \mathbb{Q}_p t_\pi^k \to \mathrm{Fil}^k B_{\mathrm{crys}}^+ \xrightarrow{\pi^{-k}F-1} B_{\mathrm{crys}}^+ \to 0$$

Note that $\mathbb{Q}_p t_\pi^k = \mathbb{Q}_p(\kappa^k)$ as $G_L$-representations. Thus the cohomology coboundary map from this exact sequence gives the reciprocity law

(1) $$L = (B_{\mathrm{crys}}^+)^{G_L} \longrightarrow H^1(L, \mathbb{Q}_p(\kappa^k))$$

Let $G = \mathrm{Gal}(L_\infty/L)$. By inflation-restriction map and local class field theory we also have the isomorphism

$$H^1(L, \mathbb{Q}_p(\kappa^k)) \cong H^1(L_\infty, \mathbb{Q}_p(\kappa^k))^G$$
$$\cong \mathrm{Hom}_G(G_{L_\infty}^{\mathrm{ab}}, \mathbb{Q}_p(\kappa^k))$$
$$\cong \mathrm{Hom}_G(U_\infty^1, \mathbb{Q}_p(\kappa^k)),$$

where $G_{L_\infty}^{\mathrm{ab}}$ is the largest abelian quotient of $G_{L_\infty}$.

In [**Ha**], Harrison proved the following explicit description of the reciprocity law from the equation (1).

**Theorem 1.** *Let $L$ be a finite unramified extension of $\mathbb{Q}_p$ and $r \geq 1$. Then, identifying $\mathbb{Q}_p(\kappa^k)$ with $\mathbb{Q}_p t_\pi^k \subseteq B_{\mathrm{crys}}^+$, the reciprocity law map from (1)*

$$\partial^k : L = (B_{\mathrm{crys}}^+)^{G_L} \to H^1(L, \mathbb{Q}_p(\kappa^k)) \cong \mathrm{Hom}(U_\infty^1, \mathbb{Q}_p(\kappa^k))$$

*is given by*

$$\partial^k(a) = \{u \mapsto \frac{1}{(k-1)!} \mathrm{Tr}_{\mathbb{Q}_p}^L (a\phi_{CW}^k(u)) t_\pi^k\}, u \in L$$

*where $\mathrm{Tr}_{\mathbb{Q}_p}^L$ is the trace map from $L$ to $\mathbb{Q}_p$.*

We next define
$$\widetilde{\log g_u}(X) = \log g_u(X) - \frac{1}{p}\sum_{\lambda \in E_\pi} \log g_u(X[+]\lambda)$$
$$= \log g_u(X) - \frac{1}{p}\log g_u^F([\pi]_\mathbf{F}(X))$$

Then $\widetilde{\log g_u}(\theta(X)) \in \mathcal{O}_{D_p}[[X]]$ corresponds to an $\mathcal{O}_{D_p}$-valued measure on $\mathbb{Z}_p$ which is actually supported on $\mathbb{Z}_p^\times \cong \mathrm{Gal}(L_\infty/L)$. Thus we can define an $\mathcal{O}_{D_p}$-valued measure $\mu_u \in \mathcal{O}_{D_p}[[G]]$ on $G = \mathrm{Gal}(L_\infty/L)$ by

$$\widetilde{\log g_u}(\theta(X)) = \int_G (1+X)^{\kappa(\sigma)} d\mu_u(\sigma)$$

and extend it to a measure on $\mathcal{G} = \mathrm{Gal}(L_\infty/\mathbb{Q}_p)$.

Following Harrison we define

$$h_u^k(X) = \left(\frac{1}{\lambda'_\mathbf{F}(X)}\frac{d}{dX}\right)^k \log g_u(X) \in \mathcal{O}_L[[X]]$$

and

$$\widetilde{h_u}^k(X) = \left(\frac{1}{\lambda'_\mathbf{F}(X)}\frac{d}{dX}\right)^k \widetilde{\log g_u}(X) \in \mathcal{O}_L[[X]]$$

The following equations are proved in [**Ha**]. For $k \geq 1$,

(2) $$\widetilde{h_u}^k(X) = h_u^k(X) - \frac{\pi^k}{p}(h_u^k)^F([\pi]_\mathbf{F}(X))$$

(3) $$\phi_{CW}^k(u) - \frac{\pi^k}{p}[\phi_{CW}^k(u)]^F = \Omega_p^{-k}\int_G \kappa(\sigma)^k d\mu_u(\sigma)$$

(4) $$\pi^{-k}(h_u^k)^{F^{-1}}(\pi_1) = (p\Omega_p)^{-k}\int_G \kappa(\sigma)^k \zeta_1^{\kappa(\sigma)} d\mu_u^{F^{-1}}(\sigma) + p^{-1}\phi_{CW}^k(u)$$

Here $\zeta_1$ is a primitive $p$-th root of unity.

We next summarize basic properties of anomality.

**Proposition 1.1.** *Let $L$ be the unramified extension of $\mathbb{Q}_p$ of degree $d \mid p-1$. Let $F$ be the Frobenius element of $\mathrm{Gal}(L/\mathbb{Q}_p)$. Let $\epsilon$ be the cyclotomic character of order $d$. For any integer $N > 0$, the following statements are equivalent.*

1. $\mu_{p^N} \subseteq L(\mathbf{F}_{\pi^\infty})$.
2. $(\frac{\pi}{p})^{[L:\mathbb{Q}_p]} \equiv 1 \bmod p^N$.
3. $\frac{\pi}{p} \equiv \epsilon(F)^S \bmod p^N$, for some $\geq 0$, where $F$ is the Frobenius element of $\mathrm{Gal}(\mathbb{Q}_p(\mu_d)/\mathbb{Q}_p)$.
   *If one (hence all) of the above conditions holds then we also have*
4. $\chi_N(\sigma) \equiv \kappa(\sigma) \bmod p^N$, for $\sigma \in \mathrm{Gal}(L(\mathbf{F}_{\pi^\infty})/L)$, where $\chi_N : \mathcal{G} \to (\mathbb{Z}/p^N\mathbb{Z})^\times$ is the character giving the action of $\mathrm{Gal}(L(\mathbf{F}_{\pi^\infty})/L)$ on $\mu_{p^N}$.
5. $\chi_N(\sigma) \equiv \epsilon^S(\sigma) \bmod p^N$, for $\sigma \in \mathrm{Gal}(L(\mathbf{F}_{\pi^\infty})/L)$ and for the integer $S$ in statement 3.

We define the index of *anomality* of $\mathbf{F}_\pi$ over $L$ to be the largest integer $N$ in the proposition. By local class field theory, this $N$ exists for $\mathbf{F}_\pi \neq \mathbb{G}_m$.

**Proof of Proposition:** (1) $\iff$ (2) is a consequence of local class field theory[**de**, I.3.7].

(2) $\iff$ (3) follows from the fact that $\mu_d$ is a subset of $\mu_{p-1} \subseteq \mathbb{Z}_p^\times$ and hence are distinct modulo $p$.

(1) $\implies$ (4) is obtained from lifting $\sigma$ to $\tilde{\sigma} \in \mathrm{Gal}(\bar{\mathbb{Q}}_p/\mathbb{Q}_p^{\mathrm{un}})$ and applying it to the equation $\pi_N = \theta^{F^{-N}}(\zeta_N - 1)$. Similarly for (5). For more details, see [**de**, **Ha**].

$\Xi$

We have the fundamental exact sequence [**de**, I.3.7] of $\mathcal{O}_{\mathrm{D}_p}[[\mathrm{Gal}(L(\mathbf{F}_{\pi^\infty})/L)]]$-modules

$$(5) \quad 0 \to U_\infty^1 \hat{\otimes}_{\mathbb{Z}_p} \mathcal{O}_{\mathrm{D}_p} \xrightarrow{i} \mathcal{O}_{\mathrm{D}_p}[[\mathrm{Gal}(L(\mathbf{F}_{\pi^\infty})/L)]] \xrightarrow{j} (\mathcal{O}_{\mathrm{D}_p}/p^N)(\chi_N) \to 0$$

where $i$ obtained by first defining $i(u) = \mu_u$ for $u \in U_\infty^1$ and then extending linearly, $j$ is given by $j(\mu) = \int_\mathcal{G} \chi_N(\sigma) d\mu(\sigma)$ regarding $\mathcal{O}_{\mathrm{D}_p}[[\mathrm{Gal}(L(\mathbf{F}_{\pi^\infty})/L)]]$ as $\mathcal{O}_{\mathrm{D}_p}$-valued measures on $\mathrm{Gal}(L(\mathbf{F}_{\pi^\infty})/L)$.

## 2. Explicit exponential maps at ordinary primes

Now we use the explicit reciprocity law (Theorem 1) to get an explicit exponential map for the Hecke character $\varphi = \psi^k \bar{\psi}^{-j}$ at a good, ordinary prime $p$, that is, at a prime $p$ where the associated elliptic curve $E/\mathbb{Q}$ has good, ordinary reduction. Fix an embedding $i_p : \bar{\mathbb{Q}} \hookrightarrow \bar{\mathbb{Q}}_p$ and let $v_p$ be the associated place of $\bar{\mathbb{Q}}$ above $p$. For $u \in \mathbb{Z}_p^\times$, let $<u>$ be the projection of $u$ onto the direct factor $1 + \mathbb{Z}_p$ of $\mathbb{Z}_p^\times$.

**Theorem 2.** *Let $p$ be an ordinary prime, $p \nmid k$. Let $\wp$ be the prime of $K$ with $v_p|\wp$ and let $\varphi_\wp$ be the $\wp$-adic character associated to $\varphi$. Let $\pi = \psi(\wp) < \bar{\psi}^{-j}(\wp) >^{1/k}$. Let $L \supseteq \mathbb{Q}_p(E_{\bar{\psi}(\wp)}^{\otimes(-j)})$ be an unramified extension. The exponential map*

$$\exp_L : L = DR_L(\mathbb{Q}_p(\varphi_\wp)) \to H^1(L, \mathbb{Q}_p(\varphi_\wp)) \cong \mathrm{Hom}_G(U_\infty^1, \mathbb{Q}_p(\varphi_\wp))$$

*of Bloch and Kato [**B-K**, §3], whose construction will be reviewed in the proof, has the following explicit description. Denote $f_a = \exp_L(a) \in \mathrm{Hom}_G(U_\infty^1, \mathbb{Q}_p(\varphi_\wp))$ with $a \in L = DR_L(\mathbb{Q}_p(\varphi_\wp))$, then*

$$f_a(u) = \frac{1}{(k-1)!} \mathrm{Tr}_{\mathbb{Q}_p}^L [(a - \frac{a^F}{\pi^k}) \phi_{CW}^k] T_\varphi.$$

*Here $\phi_{CW}^k$ is the $k$-th Coates-Wiles homomorphism (1) associated to the Lubin-Tate formal $\mathbb{Z}_p$-module for the uniformiser $\pi$, $T_\varphi$ is a $\mathbb{Q}_p$-basis of $\mathbb{Q}_p(\varphi_\wp)$ and $L$ is identified with $DR_L(\mathbb{Q}_p(\varphi_\wp))$ by $a \mapsto a t_\pi^{-k} \otimes T_k$.*

The proof of the theorem is based on the following relation between the $p$-adic Hecke characters and $p$-adic characters from Lubin-Tate formal groups.

**Proposition 2.1.** *Let $\pi$ be a uniformiser of $\mathbb{Z}_p$. Write $\pi = up$ with $u \in \mathbb{Z}_p^\times$. Let $F_\pi$ be the Lubin-Tate formal group associated to $\pi$, $\kappa_\pi : G_{\mathbb{Q}_p} \to \mathbb{Z}_p^\times$ the character obtained from the Galois action on the $\pi$-power division points in $F_\pi(\mathcal{P}_{\bar{\mathbb{Q}}_p})$. Then $\kappa_\pi = \psi_\wp \psi_{\wp^*}^r$, where $r = \frac{\log_p u}{\log_p(p/\pi_E)} + 1$, $\pi_E$ is the uniformiser for the Lubin-Tate formal group from the CM elliptic curve $E$.*

**Proof:** Since $F_\pi$ and $\hat{E} = F_{\pi_E}$ are isomorphic over the maximal unramified extension of $\mathbb{Q}_p$, $\kappa_\pi \psi_\wp^{-1}$ is an unramified character $G_{\mathbb{Q}_p} \to \mathbb{Z}_p^\times$, hence induces a character $\mathrm{Gal}(\mathbb{Q}_p^{un}/\mathbb{Q}_p) \cong \hat{\mathbb{Z}} \to \mathbb{Z}_p^\times$. Since $\psi_{\wp^*} : G_{\mathbb{Q}_p} \to \mathbb{Z}_p^\times$ is unramified and induces surjective map $\psi_{\wp^*} : \mathrm{Gal}(\mathbb{Q}_p^{nr}/\mathbb{Q}_p) \cong \hat{\mathbb{Z}} \to \mathbb{Z}_p^\times$, there is $f : \mathbb{Z}_p^\times \to \mathbb{Z}_p^\times$ such that $f \circ \psi_{\wp^*} = \kappa_\pi \psi_\wp^{-1}$. Since $\mathbb{Z}_p^\times$ is pro-cyclic, $f$ is of the form $f : x \mapsto x^r$ for some $r \in \mathbb{Z}_p^\times$. Thus $\kappa_\pi = \psi_\wp \psi_{\wp^*}^r$. By local class field theory,
$$\kappa_\pi(\pi) = \psi_\wp(\pi_E) = \psi_{\wp^*}(u) = 1.$$
By the theory of complex multiplication, $\psi_{\wp^*}(\pi_E) = p/\pi_E$. Thus $\psi_{\wp^*}(\pi) = \psi_{\wp^*}(\pi_E) = p/\pi_E$ and $\psi_\wp(\pi) = \psi_{\wp^*}(\pi)^{-r} = (p/\pi_E)^{-r}$.

On the other hand, Since $\psi_\wp(\pi_E) = 1$, we have $\psi_\wp(\pi) = \psi_\wp(\pi/\pi_E)$. Since $\hat{E}$ and $\mathbb{G}_m$ are isomorphic over the maximal unramified extension of $\mathbb{Q}_p$, $\psi_\wp \chi_p^{-1}$ is unramified. Thus we have $\psi_\wp(\pi/\pi_E) = \chi_p(\pi/\pi_E)$. By class field theory[**Ne**, p.63],
$$\xi_{p^n}^{\chi_p(\pi/\pi_E)} = (\pi/\pi_E, \mathbb{Q}_p)\xi_{p^n} = [(\pi/\pi_E)^{-1}]_{\mathbb{G}_m}\xi_{p^n} = \xi_{p^n}^{(\pi/\pi_E)^{-1}}.$$
So $\chi_p(\pi/\pi_E) = (\pi/\pi_E)^{-1}$. Thus $\frac{up}{\pi_E} = \frac{\pi}{\pi_E} = \psi_\wp(\pi) = (\frac{p}{\pi_E})^r = (\frac{p}{\pi_E})^r$ and $u = (p/\pi_E)^{r-1}$. Taking the base $p$ logarithm proves the proposition. $\Xi$

**Proof of the Theorem:** We first relate $\varphi_\wp$ to a formal group character. Let $\pi = \psi(\wp) < \bar\psi^{-j} >^{1/k}$, as defined in the theorem. Let $\kappa_\pi$ be the corresponding formal group character, defined by the Galois action of $G_{\mathbb{Q}_p}$ on the division points of the formal group $\mathbf{F}_\pi$. Since $p = \psi(\wp)\bar\psi(\wp)$ and $\psi(\wp) = \pi_E$, we have $\bar\psi(\wp) = p/\pi_E$. By Proposition 2.1, if $\pi = p(p/\pi_E)^r$ with $r \in \mathbb{Z}_p^\times$, then $\kappa_\pi = \psi_\wp \psi_{\wp^*}^r$. Let $r' \in \mathbb{Z}_p^\times$ be chosen with $(\frac{p}{\pi_E})^{r'} = <\frac{p}{\pi_E}>$, then we get
$$\pi = \psi(\wp) < \bar\psi^{-j}(\wp) >^{1/k} = p(\frac{\pi_E}{p})(\frac{\pi_E}{p})^{r'(-j/k)}.$$
Thus $\kappa_\pi = \psi_\wp \psi_{\wp^*}^{r'(-j/k)}$ and $\kappa_\pi^k = \psi_\wp^k (\psi_{\wp^*}^{r'})^{(-j)}$. Note that multiplication by $r'$ is the projection of $\mathbb{Z}_p^\times$ onto the factor $1 + \mathbb{Z}_p$. Thus $\psi_{\wp^*}^{r'}$ is trivial on the prime to $p$ part of the Galois group $\mathrm{Gal}(\mathbb{Q}_p(E_{(\wp^*)^\infty})/\mathbb{Q}_p)$ and is the identity on the pro-$p$-part of this group. Therefore $\phi_{\wp^*}^{r'}$ is the same as $\psi_{\wp^*}$ regarded as characters on $\mathbb{Q}_p(E_{\wp^*})$. Thus $\kappa^k = \psi_\wp^k \psi_{\wp^*}^{-j}$ on $L$. It follows that $t_\pi^k \otimes T_\varphi$ is an element in $DR_L(\mathbb{Q}_p(\varphi_\wp))$.

By definition, the exponential map $\mathrm{Exp}_L$ of Bloch and Kato[**B-K**] is defined to be the connecting homomorphism when we take $G_L$-invariants of the exact sequence from [**B-K**]
$$0 \to \mathbb{Q}_p(\varphi_\wp) \xrightarrow{\alpha} B_{crys}^{F=1}(\varphi_\wp) \oplus B_{DR}^+(\varphi_\wp) \xrightarrow{\beta} B_{DR}(\varphi_\wp) \to 0$$
where $\alpha$ is given by $\alpha(xT_\varphi) = (x \otimes T_\varphi, x \otimes T_\varphi)$ and $\beta$ is given by $\beta(x \otimes T_\varphi, t \otimes T_\varphi) = ((x-y) \otimes T_\varphi)$. We now choose $x \in B_{crys}^{F=1}(\varphi_\wp)$ and $y \in B_{DR}^+$ such that $x - y = at_\pi^{-k}$. Then $y = x - at_\pi^{-k} \in B_{crys}$.

Let $x_1 = t_\pi^k x$, $y_1 = t_\pi^k y$. Since $y \in Fil^0 B_{crys} \subseteq B_{crys}$ and $t_\pi^k \in B_{crys}^+$, we have
$$y_1 \in B_{crys} \cap Fil^k B_{DR} = Fil^k B_{crys}.$$
By definition, $x_1^F = \pi^k x_1$. From $x = y + at_\pi^{-k}$ we get $x_1 = y_1 + a \in Fil^k B_{crys} + B_{crys}^+ \subseteq Fil^0 B_{crys}$. Also,
$$\begin{aligned} y_1^{F^m} &= x_1^{F^m} - a^{F^m} = \pi^{km} x_1 - a^{F^m} \\ &\in Fil^{km} B_{crys} + B_{crys}^+ \subseteq Fil^0 B_{crys}, \ \forall m \geq 0.\end{aligned}$$

Now we show that $y_1 \in B_{crys}^+$. This would give $y_1 \in Fil^k B_{crys}^+$. As $B_{crys} = \cup_{s \geq 0} t^{-s} B_{crys}^+$, there is $s \geq 0$ with $y_1 \in t^{-s} B_{crys}^+$. If $s = 0$ we are done. If $s > 0$, then $t^s y_1 \in B_{crys}^+$ and

$$(t^s y_1)^{F^m} = p^{ms} t^s y_1^{F^m} \in Fil^s B_{crys} \subseteq Fil^1 B_{crys}.$$

By [**Fo1**, Prop. 4.14], this shows $t^s y_1 \in t B_{crys}^+$. Hence $t^{s-1} y_1 \in B_{crys}^+$. Repeat this process $s$ times, we get $y_1 \in B_{crys}^+$, as needed.

By definition, for $g \in G_L$,

$$\begin{aligned} \mathrm{Exp}_L(at_\pi^{-k} \otimes T_\varphi)(g) &= (x \otimes T_\varphi, y \otimes T_\varphi)^g - (x \otimes T_\varphi, y \otimes T_\varphi) \\ &= ((x \otimes T_\varphi)^g - x \otimes T_\varphi, (y \otimes T_\varphi)^g - y \otimes T_\varphi). \end{aligned}$$

which is mapped to zero via $\beta$. Hence it is an image of $\alpha$. Thus $(x \otimes T_\varphi)^g - x \otimes T_\varphi = (y \otimes T_\varphi)^g - y \otimes T_\varphi \in \mathbb{Q}_p(\varphi_\wp)$ and

$$\begin{aligned} \mathrm{Exp}_L(at_\pi^{-k} \otimes T_\varphi)(g) &= (y \otimes T_\varphi)^g - y \otimes T_\varphi \\ &= \varphi_\wp^{-1}(g) t_\pi^{-k} y_1^g \varphi_\wp \otimes T_\varphi - t_\pi^{-k} y_1 \otimes T_\varphi \\ &= t_\pi^{-k}(y_1^g - y_1) \otimes T_\varphi \end{aligned}$$

It is shown above that $y_1 \in Fil^k B_{crys}^+$, hence $(\pi^{-k} F - 1) y_1 \in B_{crys}^+$ in the exact sequence

$$0 \to \mathbb{Q}_p t_\pi^k \to Fil^k B_{\mathrm{crys}}^+ \xrightarrow{\pi^{-k} F - 1} B_{\mathrm{crys}}^+ \to 0.$$

So $y_1^g - y_1$ is the image of $(\pi^{-k} F - 1) y_1$ under the boundary map $\partial_\pi^k$. But

$$\pi^{-k} F y_1 = \pi^{-k} F(x_1 - a) = \pi^{-k} x_1^F - \pi^{-k} a^F = x_1 - \pi^{-k} a^F.$$

So

$$\begin{aligned} (\pi^{-k} F - 1) y_1 &= \pi^{-k} F(x_1 - a) - (x_1 - a) \\ &= x_1 - \pi^{-k} a^F - (x_1 - a) \\ &= a - \pi^{-k} a^F. \end{aligned}$$

Thus

$$\begin{aligned} \mathrm{Exp}_L(at_\pi^{-k} \otimes T_\varphi)(g) &= t_\pi^{-k}(y_1^g - y_1) \otimes T_\varphi \\ &= t_\pi^{-k} \partial_\pi^k((\pi^{-k} F - 1) y_1)(g) \\ &= t_\pi^{-k} \partial_\pi^k(a - \pi^{-k} a^F)(g). \end{aligned}$$

Composing it with the reciprocity map in local class field theory, and applying Theorem 1, we get, for $u \in U_\infty^1$,

$$\mathrm{Exp}_L(at_\pi^{-k} \otimes T_\varphi)(u) = \frac{1}{(k-1)!} \mathrm{Tr}_{\mathbb{Q}_p}^L[(a - \frac{a^F}{\pi^k}) \phi_{CW}^k] T_\varphi.$$

Ξ

We now consider $L/\mathbb{Q}_p$, an unramified extension of degree $d | p - 1$. For application to the Bloch-Kato conjecture, $L$ will be $\mathbb{Q}_p(E_{(\bar{\pi}_E)}^{\otimes(-j)})$. The degree is $(p-1)/gcd(p-1, j)$. By local class field theory and Kummer theory $L = \mathbb{Q}_p(\alpha)$ where $\alpha = u^{1/d}$ with $u$ a unit of $\mathbb{Q}_p^\times$ but not in $(\mathbb{Q}_p^\times)^r$ for any proper divisor $r$ of $d$. In other words, $u$ is of order $d$ in $\mathbb{Q}_p^\times/(\mathbb{Q}_p^\times)^d$. As usual, for $1 \leq n \leq \infty$, define $L_n = L(\mathbf{F}_{\pi^n})$. Define $\mathcal{G} = \mathrm{Gal}(L_\infty/\mathbb{Q})$, $G = \mathrm{Gal}(L_\infty/L)$, $\Delta = \mathrm{Gal}(L(\mathbf{F}_\pi)/L)$, $\Delta_1 = \mathrm{Gal}(L(\mathbf{F}_\pi)/\mathbb{Q}_p)$

and $\Delta' = \mathrm{Gal}(L/\mathbb{Q}_p)$. Then we have canonical isomorphisms $\mathcal{G} \cong G \times \Delta' \cong \Gamma \times \Delta_1$, $G \cong \Gamma \times \Delta$, $\Delta_1 \cong \Delta \times \Delta'$. The related fields are shown in the following diagram.

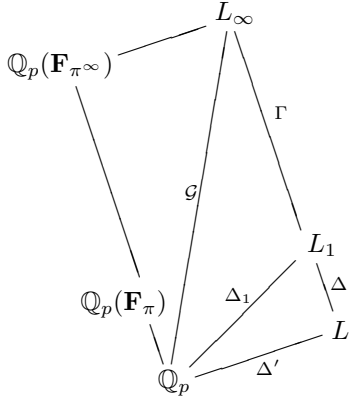

Let $\omega = \kappa_\pi|_\Delta$. Let $\epsilon$ be the character induced from the action of $\Delta'$ on $\alpha$. It is in fact the cyclotomic character of order $d$. These characters can be regarded as characters of $\mathcal{G}$ acting trivially on $\Gamma$. From our choice of $d$, $\mu_d \subseteq \mu_{p-1} \subseteq \mathbb{Z}_p^\times$ and are distinct modulo $p$. It follows from Proposition 1.1 that the index of anomaly $N$ of $L$ is the largest $N$ such that $\pi/p \equiv \epsilon(F)^S \mod p^N$ for some $S$, where $F$ is the Frobenius of $\Delta'$. If $N > 0$ we choose $S$ such that $0 \leq S \leq d - 1$.

We have canonical decomposition $U_\infty^1 \cong \oplus_{(s,t) \bmod (d,p-1)} (U_\infty^1)^{(s,t)}$ of $U_\infty^1$ on which $\Delta_1$ acts via $\epsilon^s \omega^t$. We can also decompose $L_1 \cong \oplus_{0 \leq t \leq p-2} L_1^{(t)}$ into the 1-dimensional $L$-eigenspaces of $L_1$ with eigen-character $\omega^t$. Choose a $L$-basis $e_t$ of $L_1^{(t)}$ such that $e_t = 1$ if $t = 0$ and $0 \leq v_{\mathfrak{P}}(e_t) \leq p - 2$ for any $t$, where $\mathfrak{P}$ is the maximal ideal of $\mathcal{O}_{L_1}$. This could be done by changing $e_t$ by a power of $p$ if necessary. Now $\mathbf{D}_p(\mathbf{F}_\pi) = \mathbf{D}_p(\mu_p)$ and if we identify $\mathrm{Gal}(\mathbf{D}_p(\mathbf{F}_\pi)/\mathbf{D}_p)$ with $\mathrm{Gal}(L(\mathbf{F}_\pi)/L) = \Delta$ then $\omega$ is the cyclotomic character from the action on $\mu_p$. Thus we can define the Gauss sum

$$G(t) = \begin{cases} \frac{1}{p-1} \sum_{\sigma \in \Delta} \omega^{-t}(\sigma) \zeta_1^\sigma, & t \not\equiv 0 \\ 1, & t = 0 \end{cases}$$

It is well-known that $v_{\mathfrak{P}_1}(G(t)) = t$, $0 \leq t \leq p - 2$, where $\mathfrak{P}_1$ is the maximal ideal of $\mathbf{D}_p(\mathbf{F}_\pi)$. Since $G(t)$ and $e_t$ are both in $\mathbf{D}_p(\mathbf{F}_\pi)^{(t)}$ and have $\mathfrak{P}_1$-adic valuation between $0$ and $p - 2$, it follows that they must differ by a unit: $G(t) = u_t e_t$, $u_t \in \mathcal{O}_{\mathbf{D}_p}^\times$, $0 \leq t \leq p - 2$. As $L \cong \oplus_{0 \leq s \leq d-1} \mathbb{Q}_p \alpha^s$, we get the decomposition $L_1 \cong \oplus_{(s,t) \bmod (d,p-1)} L_1^{(s,t)}$ of $L_1$ into eigenspaces $L_1^{(s,t)} = \mathbb{Q}_p \alpha^s e_t$ of $L_1$ with eigen-character $\epsilon^s \omega^t$.

For $k \geq 1$, define

$$\phi_k : U_\infty^1 \to L_1, \ u \mapsto \pi^{-k} (h_u^k)^{F^{-1}}(\pi_1).$$

Recall from section 1 that $\pi_1 = \theta^{F^{-1}}(\zeta_1 - 1)$, and $\theta(X) \in \mathcal{O}_{\mathbf{D}_p}[[T]]$ is an $\mathbb{Z}_p$-isomorphism from $\mathbb{G}_m$ onto $\mathbf{F}_\pi$. $\phi_k$ is a $\mathbb{Z}_p$-homomorphism. We will explicitly determine its image. This will be play an essential role in application to the Bloch-Kato conjecture.

**Proposition 2.2.** *Let $k \geq 2$. Define $N_k = \min\{N, 1 + v_p(k-1)\}$. Then*
$$\phi_k(U_\infty^1) = p^{-k}\left[\oplus_{s \neq -S \text{ or } t \not\equiv 1-k} \mathbb{Z}_p \alpha^s e_t \oplus p^{N_k}\mathbb{Z}_p \alpha^{-S} e_{1-k}\right]$$
*If $\mathcal{L}_k$ is the lattice of $L_1$ which is the exact annihilator to $\frac{1}{(k-1)!}\phi_k(U_\infty^1)$ under the pairing*
$$Tr_{\mathbb{Q}_p}^{L_1}: L_1 \times L_1 \to \mathbb{Q}_p/\mathbb{Z}_p$$
*then*
$$\mathcal{L}_k = \begin{cases} p^k(k-1)!\left[\bigoplus_{t \equiv 0(p-1)} \mathbb{Z}_p \alpha^s e_t \oplus \bigoplus_{t \not\equiv 0(p-1)} p^{-1}\mathbb{Z}_p \alpha^s e_t\right], & \text{if } N = 0 \\ p^k(k-1)!\left[\bigoplus_{t \equiv 0(p-1)} \mathbb{Z}_p \alpha^s e_t \right. \\ \quad \oplus \bigoplus_{t \not\equiv 0, k-1(p-1), \text{or } t \equiv k-1, s \neq S} p^{-1}\mathbb{Z}_p \alpha^s e_t \\ \left.\oplus p^{-N_k-1}\mathbb{Z}_p \alpha^S e_{k-1}\right], & \text{if } N \neq 0 \end{cases}$$

**Proof**: This is again similar to the $j = 0$ case. Recall that
$$h_u^k(X) = \left(\frac{1}{\lambda_{\mathbf{F}}'(X)}\frac{d}{dX}\right)^k \log g_u(X) \in \mathcal{O}_L[[X]].$$

As in [**de**, I.3.5.3], the properties of $g_u$ imply that $\phi_k(u^\sigma) = \kappa(\sigma)^k[\phi_k(u)]^\sigma$ for $\sigma \in \mathcal{G}$ where we consider $\kappa$ as a character on $\mathcal{G}$ that is trivial on $\Delta'$. Therefore if $\sigma \in \Delta$, we have $\phi_k(u^\sigma) = \omega^k(\sigma)[\phi_k(u)]^\sigma$ and so $\phi_k((U_\infty^1)^{(s,t)})$ is contained in $L_1^{(s,t-k)} = \mathbb{Q}_p \alpha^s e_{t-k}$.

As $\phi_k$ is $\mathbb{Z}_p$-linear we have $\phi_k((U_\infty^1)^{(s,t)}) = p^{r_{s,t}}\mathbb{Z}_p \alpha^s e_{t-k}$ for some $r_{s,t} \in \mathbb{Z}$ and so by $G(t) = u_t e_t$, $u_t \in \mathcal{O}_{\mathbf{D}_p}^\times$
$$\phi_k((U_\infty^1)^{(s,t)})\mathcal{O}_{\mathbf{D}_p} = p^{r_{s,t}}G(t-k)\mathcal{O}_{\mathbf{D}_p}$$
as a subset of $\mathbf{D}_p(\mathbf{F}_\pi)$. We only need to compute the integers $r_{s,t}$.

We first decompose $\mathcal{O}_{\mathbf{D}_p}[[G]]$ into eigenspaces for the eigen-characters $\omega^t$ of $\Delta$
$$\mathcal{O}_{\mathbf{D}_p}[[G]] = \oplus_{0 \leq t \leq p-2}\mathcal{O}_{\mathbf{D}_p}[[\Gamma]]E_t,$$
with $E_t = \frac{1}{p-1}\sum_{\sigma \in \Delta}\omega^{-t}(\sigma)\sigma \in \mathbb{Z}_p[\Delta]$. This enable us to write $\mu_u$ as $\mu_u = \sum E_t\mu_{u,t}$, $\mu_{u,t} \in \mathcal{O}_{\mathbf{D}_p}[[\Gamma]]$, considered as elements in $\mathcal{O}_{\mathbf{D}_p}[[G]]$. For each $0 \leq t \leq p-2$, define
$$\psi_t: U_\infty^1 \to \mathcal{O}_{\mathbf{D}_p}, \quad u \mapsto \int_\Gamma \kappa(\sigma)^k d\mu_{u,t}(\sigma).$$

This map is $\mathbb{Z}_p$-linear and hence can be extended to a $\mathcal{O}_{\mathbf{D}_p}$-linear map from $U_\infty^1 \hat{\otimes}_{\mathbb{Z}_p} \mathcal{O}_{\mathbf{D}_p}$ to $\mathcal{O}_{\mathbf{D}_p}$. The decomposition of $U_\infty^1$ gives a decomposition of $U_\infty^1 \hat{\otimes}_{\mathbb{Z}_p} \mathcal{O}_{\mathbf{D}_p}$ into the sum of eigenspaces $(U_\infty^1)^{(s,t)} \hat{\otimes}_{\mathbb{Z}_p} \mathcal{O}_{\mathbf{D}_p}$. From $\mu_{u^\sigma} = \sigma\mu_u$, $\sigma \in G$, we see that for $\sigma \in \Delta$, $\mu_{u^\sigma,t} = \omega^t(\sigma)\mu_{u,t}$ and $\psi_t(u^\sigma) = \omega^t(\sigma)\psi_t(u)$. Thus for $u \in (U_\infty^1)^{(s,t)} \hat{\otimes}_{\mathbb{Z}_p} \mathcal{O}_{\mathbf{D}_p}$, we have $\psi_{t_1}(u) = 0$ unless $t_1 = t$. Similarly from $\phi_{CW}^k(u^\sigma) = \kappa^k(\sigma)\phi_{CW}^k(u) = \omega^k(\sigma)\psi_{CW}^k(u)$, $\sigma \in \Delta$, we obtain that $\phi_{CW}^k(u) = 0$ for $u \in (U_\infty^1)^{(s,t)} \hat{\otimes}_{\mathbb{Z}_p} \mathcal{O}_{\mathbf{D}_p}$ unless $k \equiv t \mod p-1$.

By equation (4),
$$\phi_k(u) = (p\Omega_p)^{-k}\int_G \kappa(\sigma)^k \zeta_1^{\kappa(\sigma)} d\mu_u^{F^{-1}}(\sigma) + p^{-1}\phi_{CW}^k(u)$$
$$= (p\Omega_p)^{-k}\frac{1}{p-1}\sum_t \sum_{\sigma \in \Delta}\omega^{-t}(\sigma)\int_\Gamma \kappa(\sigma\gamma)^k \zeta_1^{\kappa(\sigma\gamma)} d\mu_{u,t}^{F^{-1}}(\gamma) + p^{-1}\phi_{CW}^k(u)$$

Since $\kappa(\sigma\gamma) = \kappa(\sigma)\kappa(\gamma) = \omega(\sigma)\kappa(\gamma)$, and for $\gamma \in \Gamma$, we have $\kappa(\gamma) \equiv 1 \mod p$ and $\zeta_1^{\kappa(\sigma\gamma)} = \zeta_1^{\omega(\sigma)}$, we obtain

$$\phi_k(u) = (p\Omega_p)^{-k} \sum_t G(t-k)[\psi_t(u)]^{F^{-1}} + p^{-1}\phi_{CW}^r(u) \tag{6}$$

Now if $u \in (U_\infty^1)^{(s,t)}\hat{\otimes}_{\mathbb{Z}_p}\mathcal{O}_{\mathbf{D}_p}$ with $t \not\equiv k \mod p-1$, then $\phi_{CW}^k(u) = 0$ and $\psi_{t_1}(u) = 0$ for $t_1 \neq t$. So the previous equation implies that $\phi_k(u) = (p\Omega_p)^{-k}G(t-k)[\psi_t(u)]^{F^{-1}}$. Hence

$$\phi_k((U_\infty^1)^{(s,t)})\mathcal{O}_{\mathbf{D}_p} = p^{-k}G(t-k)[\psi_t((U_\infty^1)^{(s,t)}\hat{\otimes}_{\mathbb{Z}_p}\mathcal{O}_{\mathbf{D}_p})]^{F^{-1}}.$$

On the other hand, if $u \in (U_\infty^1)^{(s,t)}\hat{\otimes}_{\mathbb{Z}_p}\mathcal{O}_{\mathbf{D}_p}$ with $t \equiv k \mod p-1$, then it is easy to see that $\phi_{CW}^k(u^\tau) = [\phi_{CW}^k(u)]^\tau$ for $\tau \in \Delta' = G(L/\mathbb{Q}_p) \cong G(L(\mathbb{F}_{\pi^\infty})/\mathbb{Q}_p(\mathbb{F}_{\pi^\infty})))$, so $\phi_{CW}^k(u) \in L^s$, the eigenspace for the character $\chi^s$, and

$$[\phi_{CW}^k(u)]^F = \chi^s(F)\phi_{CW}^k(u) = \rho^s\phi_{CW}^k(u).$$

Then equation (3) shows that

$$(1 - \frac{\pi^k}{p}\rho^s)\phi_{CW}^k(u) = \Omega_p^{-k}\int_G \kappa(\sigma)^k d\mu_u(\sigma) = \Omega_p^{-k}\psi_r(u) = \Omega_p^{-k}\psi_t(u).$$

Then equation (6) gives

$$\phi_k(u) = (p\Omega_p)^{-k}[\psi_t(u)]^{F^{-1}} + p^{-1}\Omega_p^{-k}(1 - \frac{\pi^k}{p})^{-1}\psi_t(u).$$

It remains to determine $\psi_t((U_\infty^1)^{(s,t)}\hat{\otimes}_{\mathbb{Z}_p}\mathcal{O}_{\mathbf{D}_p})$. Since $\chi_N(\sigma) \equiv \epsilon^{-S}(\sigma) \mod p^N$ by Proposition 1.1, the exact sequence (5)

$$0 \to U_\infty^1\hat{\otimes}_{\mathbb{Z}_p}\mathcal{O}_{\mathbf{D}_p} \xrightarrow{i} \mathcal{O}_{\mathbf{D}_p}[[\mathcal{G}]] \xrightarrow{j} (\mathcal{O}_{\mathbf{D}_p}/p^N)(\chi_N) \to 0$$

induces exact sequences of $\mathcal{O}_{\mathbf{D}_p}[[\Gamma]]$-modules

$$0 \to (U_\infty^1)^{(s,t)}\hat{\otimes}_{\mathbb{Z}_p}\mathcal{O}_{\mathbf{D}_p} \to \mathcal{O}_{\mathbf{D}_p}[[\Gamma]] \to 0$$

if $(s,t) \not\equiv (-S,1) \mod (d, p-1)$ and

$$0 \to (U_\infty^1)^{(s,t)}\hat{\otimes}_{\mathbb{Z}_p}\mathcal{O}_{\mathbf{D}_p} \to \mathcal{O}_{\mathbf{D}_p}[[\Gamma]] \to (\mathcal{O}_{\mathbf{D}_p}/p^N)(\kappa) \to 0$$

if $(s,t) \equiv (-S,1) \mod (d, p-1)$. If $(s,t) \not\equiv (-S,1) \mod (d,p-1)$, then

$$\psi_t((U_\infty^1)^{(s,t)}\hat{\otimes}_{\mathbb{Z}_p}\mathcal{O}_{\mathbf{D}_p}) = \mathcal{O}_{\mathbf{D}_p}.$$

So $\phi_k((U_\infty^1)^{(s,t)})\mathcal{O}_{\mathbf{D}_p} = p^{-k}G(t-k)\mathcal{O}_{\mathbf{D}_p}$ which means that

$$\phi_k((U_\infty^1)^{(s,t)}) = p^{-k}\mathbb{Z}_p\alpha^s e_{t-k}.$$

If $(s,t) \equiv (-S,1) \mod (d,p-1)$, then the above exact sequence shows that

$$\psi_1((U_\infty^1)^{(-S,1)}\hat{\otimes}_{\mathbb{Z}_p}\mathcal{O}_{\mathbf{D}_p})\mathcal{O}_{\mathbf{D}_p}$$
$$= \{P(\kappa(\gamma_0)^k - 1) \mid P(X) \in \mathcal{O}_{\mathbf{D}_p}[[X]], P(\kappa(\gamma_0) - 1) \equiv 0 \mod p^N\}$$

when we identify $\mathcal{O}_{\mathbf{D}_p}[[\Gamma]]$ with $\mathcal{O}_{\mathbf{D}_p}[[X]]$ by choosing a generator $\gamma_0$ of $\Gamma$ with $\kappa(\gamma_0) = 1 + p$ and identifying $\gamma_0$ with $1 + X$. Let $S_n$ be the set of polynomials of degree $n$ in $\mathcal{O}_{D_p}[[X]]$ such that $P(p) \equiv 0 \mod p^N$, then it is easily shown that $P((1+p)^k - 1) \equiv 0 \mod p^{N_k}$ for any $P \in S_n$ and that there is a $P \in S_N$ such that $P((1+p)^k-1) \not\equiv 0 \mod p^{N_k+1}$. Since every element of $\mathcal{O}_{\mathbf{D}_p}[[X]]$ can be written as a product of a polynomial and a unit power series, the same is true for the set of power series whose polynomial part has degree $n$. Therefore $\psi_1((U_\infty^1)^{(-S,1)}\hat{\otimes}_{\mathbb{Z}_p}\mathcal{O}_{\mathrm{D}_p}) = p^{N_k}$. Thus $\phi_k((U_\infty^1)^{(-S,1)}) = p^{N_k-k}\mathbb{Z}_p\alpha^{-S}e_{1-k}$. This proves the first statement.

To prove the second statement, first note that $Tr^{L_1}_{\mathbb{Q}_p}(\alpha^s e_t) = 0$ unless $s = t \equiv 0$. It then follows that $Tr^{L_1}_{\mathbb{Q}_p}(\alpha^s e_t \cdot \alpha^{s_1} e_{t_1}) = 0$ unless $s \equiv -s_1$ and $t \equiv -t_1$ since $\alpha^s e_t \cdot \alpha^{s_1} e_{t_1} \in \alpha^{s+s_1} e_{t+t_1}$ by definition. Finally, $e_1 = 1$, and, for $t \neq 0$, since $v_\wp(e_t) = v_\wp(G(t)) = t$, where $\wp$ is the prime of $L_1$, we get $v_\wp(e_t \cdot e_{-t}) = p - 1$ and hence $v_p(e_t \cdot e_{-t}) = 1$. The second statement follows from these facts and the first statement. $\boxempty$

DEPARTMENT OF MATHEMATICS AND COMPUTER SCIENCE, RUTGERS UNIVERSITY, NEWARK, NJ 07102

*E-mail address*: `liguo@andromeda.rutgers.edu`